\documentclass{article}[11pt]
\usepackage{amsthm,amsmath,a4wide}
\usepackage{amssymb}
\input xypic
\usepackage[all,tips]{xy}

\newtheorem{thm}{Theorem}[section]
\newtheorem{lem}[thm]{Lemma}
\newtheorem{prop}[thm]{Proposition}
\newtheorem{cor}[thm]{Corollary}
\newtheorem{conj}[thm]{Conjecture}

\theoremstyle{definition}
\newtheorem{defin}[thm]{Definition}

\theoremstyle{remark}
\newtheorem{remark}[thm]{Remark}
\newtheorem{example}[thm]{Example}

\newcommand{\bth}{\begin{thm}}
\renewcommand{\eth}{\end{thm}}
\newcommand{\bpr}{\begin{prop}}
\newcommand{\epr}{\end{prop}}
\newcommand{\ble}{\begin{lem}}
\newcommand{\ele}{\end{lem}}
\newcommand{\bco}{\begin{cor}}
\newcommand{\eco}{\end{cor}}
\newcommand{\bde}{\begin{defin}}
\newcommand{\ede}{\end{defin}}
\newcommand{\bex}{\begin{example}}
\newcommand{\eex}{\end{example}}
\newcommand{\bre}{\begin{remark}}
\newcommand{\ere}{\end{remark}}
\newcommand{\bcj}{\begin{conj}}
\newcommand{\ecj}{\end{conj}}

\newcommand{\beq}{\begin{equation}}
\newcommand{\eeq}{\end{equation}}
\newcommand{\ve}{{\varepsilon}}
\newcommand{\ot}{{\otimes}}
\newcommand{\op}{{\oplus}}
\newcommand{\lb}{\label}
\newcommand{\nl}{\newline}
\newcommand{\bpf}{\begin{proof}}
\newcommand{\epf}{\end{proof}}

\newcommand{\C}{{\cal C}}
\newcommand{\Z}{{\cal Z}}
\newcommand{\D}{{\cal D}}
\newcommand{\bZ}{{\mathbb Z}}

\newcommand{\bF}{{\mathbb F}}
\newcommand{\V}{{\cal V}}
\newcommand{\fsl}{{\mathfrak sl}}
\newcommand{\e}{{\mathfrak e}}

\newcommand{\Rep}{{\cal R}{\it ep}}
\newcommand{\Vect}{{\cal V}{\it ect}}
\begin{document}
\author{Alexei Davydov}
\title{Bogomolov multiplier, double class-preserving automorphisms and modular invariants for orbifolds}
\maketitle
\date{}

\begin{center}
Department of Mathematics, Ohio University, Athens, OH 45701, USA
\end{center}

\begin{abstract}
We describe the group $Aut_{br}^1(\Z(G))$ of braided tensor autoequivalences of the Drinfeld centre of a finite group $G$ isomorphic to the identity functor (just as a functor). We prove that the semi-direct product $Out_{2-cl}(G)\ltimes B(G)\ $ of the group of double class preserving automorphisms and the Bogomolov multiplier of $G$ is a subgroup of $Aut_{br}^1(\Z(G))$. An automorphism of $G$ is double class preserving if it preserves conjugacy classes of pairs of commuting elements in G.
The Bogomolov multiplier $B(G)$ is the subgroup of its Schur multiplier $H^2(G,k^*)$ of classes vanishing on abelian subgroups of $G$.
We show that elements of $Aut^1_{br}(\Z(G))$ give rise to different realisations of the charge conjugation modular invariant for $G$-orbifolds of holomorphic conformal field theories.
\end{abstract}
\tableofcontents
\section{Introduction}

The aim of this note is to construct examples of different rational conformal field theories with the same chiral algebras and the same (charge conjugation) modular invariant.
\nl
The state space of a (2-dimesional) conformal field theory comes equipped with amplitudes associated to 
a finite collection of fields inserted into a surface \cite{dms,MS}. 
Fields, whose amplitudes depend (anti-)holomorphically on insertion 
points, form what is known as (anti-)chiral algebra of the CFT. The state 
space is naturally a representation of the product of the chiral and anti-chiral 
algebras.
A conformal field theory is {\em rational} if the state space  is a finite sum of tensor 
products of irreducible representations of chiral and anti-chiral algebras.
The matrix of multiplicities of irreducible representations in the decomposition of the state space is called the {\em modular invariant} of the RCFT.
The simplest case (the so-called {\em Cardy case}) is the case of the {\em charge conjugation} modular invariant, which assumes that chiral and anti-chiral algebras coincide.
\nl 
The name modular invariant comes from the fact that the matrix of multiplicities is invariant with respect to the modular group action on characters. This fact was used in \cite{ciz} to classify modular invariants for affine $sl(2)$ rational conformal field theories.
This paper started an activity aimed at classifying possible modular invariants for various conformal field theories.
It took some time to realise that not all modular invariants correspond to conformal field theories \cite{fss,frs0}.
In this paper we show that there are different rational conformal field theories with the same (charge conjugation) modular invariant.
Thus although being a convenient numerical invariants of a rational conformal field theory modular invariants are far from being complete.
\nl
An adequate description of rational conformal field theories was obtained relatively recently (see \cite{frs} and references therein).
Mathematical axiomatisation of chiral algebras in CFTs is the notion of {\em vertex operator algebra} \cite{Bo,FHL,Ka}. A vertex operator algebra is called {\em rational} if it is a chiral algebra of a RCFT, in particular its category of modules is semi-simple. In this case the category of modules has more structure (see \cite{Hu} and references therein), it is the so-called {\em modular category}. This type of tensor categories was first studied in physics \cite{MS} and then axiomatised mathematically in \cite{Tu}. 
The state space of a RCFT corresponds to a special commutative algebra in the product or categories of modules of chiral and anti-chiral algebras (see \cite{frs1,frs2,frs3,frs4,frs5,ffrs,HK} and references therein).
These special class of commutative algebras in braided tensor categories is known as {\em Lagrangian} algebras \cite{dmno}. 
The modular invariant is just the class of this Lagrangian algebra in the Grothendieck ring.

Here we give examples of Lagrangian algebras with the charge conjugation modular invariant.
It is straightforward to see that such algebras should correspond to braided tensor autoequivalences of the category of representations of one of the chiral algebra.
A braided tensor autoequivalence corresponds to the charge conjugation modular invariant if it does not permute (isomorphism classes of) simple objects. We call such braided tensor autoequivalence {\em soft}.
We provide examples of modular categories with soft braided tensor autoequivalences.

Our examples come from permutation orbifolds of {\em holomorphic} conformal field theories (CFTs whose state space is an irreducible module over the chiral algebras).
It is argued in \cite{Ki} (see also \cite{Mu}) that the modular category of the $G$-orbifold of a holomorphic conformal field theory is the so called {\em Drinfeld} (or {\em monoidal}) {\em centre} $\Z(G,\alpha)$, where $\alpha$ is a 3-cocycle of the group $G$.
It is also known that the cocycle $\alpha$ is trivial for {\em permutation orbifolds} (orbifolds where the group $G$ is a subgroup of the symmetric group permuting copies in a tensor power of a holomorphic theory).
The assumption crucial for the arguments of \cite{Ki}  is the existence of twisted sectors.
This assumption is known to be true for permutation orbifolds \cite{bdm}. 

Clearly the property of a braided tensor autoequivalence to preserve the isomorphism classes of simple objects is stable under the composition of tensor autoequivalences. 
Thus the (tensor) isomorphism classes of soft braided tensor autoequivalences of a category $\C$ form a group $Aut^1_{br}(\C)$.
In this paper we describe the group $Aut^1_{br}t(\Z(G))$ as a semi-direct product:
$$Aut^1_{br}(\Z(G))\ \simeq\ Out_{2-cl}(G)\ltimes B(G)\ .$$
The subgroup $B(G)$ is the so-called {\em Bogomolov multiplier} of a finite group $G$. It is a subgroup of the Schur multiplier $H^2(G,k^*)$ defined as the kernel of the restriction map
$$B(G) = ker\big(H^2(G,k^*) \stackrel{\small res}{\longrightarrow}  \bigoplus_{A\subset G} H^2(A,k^*)\big)$$
where the direct sum is taken over all abelian subgroups of $G$. 
\nl
The group $Out_{2-cl}(G)$ is a subgroup of the group $Out(G)$ of the outer automorphisms of $G$ (classes of automorphisms modulo inner automorphisms). The subgroup $Out_{2-cl}(G)$ consists of outer automorphisms preserving conjugacy classes of pairs of commuting elements in $G$. 
The action of $Out_{2-cl}(G)$ on $B(G)$ comes from the natural action of $Out(G)$ on $H^2(G,k^*)$. 

The Bogomolov multiplier is an important invariant of birational geometry \cite{bog}. Recent activity produced examples of finite groups with non-trivial Bogomolov multiplier (see for example \cite{chkk,jm}). 
\nl
Finite groups with non-trivial $Out_{2-cl}(G)$ were constructed in \cite{Sz}. 

The paper is organised as follows.
We start by defining soft braided tensor autoequivalences and collect some of their basic properties in section \ref{sta}.
Then we recall necessary facts about tensor autoequivalences of the categories of group-graded vector spaces and of representations of groups (section \ref{gvs}, \ref{sar}).
We proceed by proving our main result - the description of the the group of soft braided tensor autoequivalences of the Drinfeld centre $\Z(G)$ (section \ref{sbadc}). 
Examples of non-trivial soft braided tensor autoequivalences are given in section \ref{exa}.
We conclude with a sketch of a construction of non-trivial conformal field theories with charge conjugation modular invariant (section \ref{app}).

Throughout $k$ denote an algebraically closed field of characteristic zero.
We will work with tensor categories.
A category $\C$ is {\em tensor} over $k$ if it is monoidal and enriched in the category $\Vect_k$ of finite dimensional vector spaces over $k$. That is hom-sets $\C(X,Y)$ for $X,Y\in\C$ are finite dimensional vector spaces over $k$ and the composition and tensor product of morphisms are bilinear maps.
\nl
A tensor category is {\em fusion} if it is semi-simple with finitely many simple objects.
Note that for a semi-simple $\C$ the natural embedding of the Grothendieck group into the dual of the endomorphism algebra of the identity functor
$$K_0(\C)\ \to\ End(Id_\C)^*,\qquad X\mapsto (a\mapsto a_X),\qquad X\in \C,\ a\in End(Id_\C)$$
induces an isomorphism
\beq\lb{gsc}
K_0(\C)\ot_\bZ k\ \to\ End(Id_\C)
\eeq

\section*{Acknowledgment}

The author would like to thank D. Nikshych, V. Ostrik and C. Schweigert for very useful comments on the preliminary version of the paper.

\section{Group-theoretical braided fusion categories and their autoequivalences}

\subsection{Soft tensor autoequivalences}\lb{sta}

Let $\C$ be a tensor category.
A tensor autoequivalence $F:\C\to\C$ is called {\em soft} if it is isomorphic as just a $k$-linear functor to the identity functor $Id_\C$.

\bre\lb{twf}
Any soft tensor autoequivalence $F:\C\to\C$ is a {\em twisted form} of the identity functor $Id_\C$ in the terminology of \cite{da0}.
That is $F$ is tensor isomorphic to the identity functor $Id_\C$ with the tensor structure given by a collection of automorphisms $\gamma_{X,Y}:X\ot Y\to X\ot Y$ natural in $X,Y\in \C$ and satisfying the normalised 2-cocycle condition (see \cite{da0} for details). 
\ere

For a tensor category $\C$ denote by $Aut_\ot(\C)$ the group of tensor isomorphism classes of tensor autoequivalences of $\C$. 
Denote by $Aut^1_\ot(\C)$ the subgroup of $Aut_\ot(\C)$ consisting of isomorphism classes of soft tensor autoequivalences. Clearly the subgroup $Aut^1_\ot(\C)\subset Aut_\ot(\C)$ is normal, since it is the kernel of the homomorphism forgetting tensor structure
$$Aut_\ot(\C) \to Aut(\C)\ .$$
Here $Aut(\C)$ is the group of isomorphism classes of $k$-linear autoequivalences of $\C$. 

\bre\lb{fuc}
Note that for a fusion category $\C$ a tensor autoequivalence $F:\C\to\C$ is soft if and only if it induces a trivial automorphism of the Grothendieck ring $K_0(F)=1:K_0(\C)\to K_0(\C)$.
In other words $Aut^1_\ot(\C)$ is the kernel of the natural action of $Aut_\ot(\C)$ on the Grothendieck ring $K_0(\C)$:
$$Aut^1_\ot(\C) = ker\big(Aut_\ot(\C)\to Aut(K_0(\C))\big)\ .$$
Due to the isomorphism \eqref{gsc} $Aut^1_\ot(\C)$ is also the kernel of the natural action of $Aut_\ot(\C)$ on the endomorphism algebra of the identity functor
$$Aut^1_\ot(\C) = ker\big(Aut_\ot(\C)\to Aut(End(Id_\C))\big)\ .$$
\ere

Similarly for a braided tensor category $\C$ denote by $Aut_{br}(\C)$ the group of tensor isomorphism classes of braided tensor autoequivalences of $\C$. 
Denote by $Aut^1_{br}(\C)$ the normal subgroup of $Aut_{br}(\C)$ consisting of isomorphism classes of soft braided tensor autoequivalences. 

\bre\lb{bfc}
Similarly to the remark \ref{fuc} for fusion $\C$ the group $Aut^1_{br}(\C)$ is the kernel of the natural action of $Aut_{br}(\C)$ on the Grothendieck ring $K_0(\C)$:
$$Aut^1_{br}(\C) = ker\big(Aut_{br}(\C)\to Aut(K_0(\C))\big)\ $$
and on the endomorphism algebra of the identity functor:
$$Aut^1_{br}(\C) = ker\big(Aut_{br}(\C)\to Aut(End(Id_\C))\big)\ .$$
\ere

We say that braided tensor autoequivalence $F:\C\to\C$ preserves a full tensor subcategory $\D\subset\C$ if there is a braided tensor autoequivalence $F|_\D:\D\to\D$ such that the diagram of functors
$$\xymatrix{\D \ar[rr] \ar[d]_{F|_\D} && \C \ar[d]^F\\ \D \ar[rr] && \C}$$
commutes up to a braided natural isomorphism.

The following is straightforward.
\ble\lb{rsa}
Let $\D$ be a full tensor subcategory of a braided tensor category $\C$.
Then any soft braided tensor autoequivalence $F$ of $\C$ preserves $\D$. 
\ele
In particular there is defined a homomorphism of groups
\beq\lb{rse} Aut^1_{br}(\C) \to Aut^1_{br}(\D),\qquad F\mapsto F|_\D\ .\eeq

\subsection{Soft tensor autoequivalences of graded vector spaces}\lb{gvs}

Here we recall a well-known description of tensor autoequivalences for categories of graded vector spaces. 

Let $G$ be a group.
A {\em $G$-graded} vector space $V$ is a direct sum $\op_{g\in G}V_g$ of vector spaces labelled by elements of $G$.
A morphism of $G$-graded vector spaces $U\to V$ is the direct sum $\op_{g\in G}f_g$ of homomorphisms $f_g:U_g\to V_g$.
\nl
Denote by $\V(G)$ the category of finite dimensional $G$-graded vector spaces.
\nl
Define the $G$-grading on the tensor product $U\ot V$ of $G$-graded vector spaces:
$$(U\ot V)_g = \bigoplus_{g_1g_2=g}U_{g_1}\ot V_{g_2}\ .$$

Denote by $I\in \V(G)$ the one-dimensional vector space concentrated in the trivial degree
$$I_e = k,\qquad I_g = 0,\quad g\not=e\ .$$
The category $\V(G)$ is tensor.

More generally for $f\in G$ denote by $I(f)\in \V(G)$ the one-dimensional vector space concentrated in the trivial degree
$$I(f)_f = k,\qquad I_g = 0,\quad g\not=f\ .$$
Clearly $I(f)$ is a simple object of $\V(G)$. 
The category $\V(G)$ is fusion and any simple object of $\V(G)$ is isomorphic to $I(f)$ for some  $f\in G$.
Tensor product of simple objects has the form
$$I(f)\ot I(g)\simeq I(fg)\ .$$
The Grothendieck ring of $\V(G)$ is the group ring of $G$:
$$K_0(\V(G))\ \simeq\ \bZ[G]\ .$$
A tensor autoequivalence $F$ of $\V(G)$ induces an automorphism of the group of isomorphism classes of simple objects.
Hence we have a homomorphism
\beq\lb{aos}Aut_\ot(\V(G))\to Aut(G)\ .\eeq

Let $\phi:G\to F$ be a homomorphism of groups.
For an $F$-graded vector space $V=\op_{f\in F}V_f$ define the $G$-grading on $V$ by $V_g = V_{\phi(g)}$.
We call the $G$-graded vector space $\phi^*(V) = \op_{g\in G}V_g$ the {\em inverse image} of $V$ along the homomorphism $\phi:G\to F$.
The inverse image functor
$$\phi^*:\V(F)\to\V(G)\qquad V\mapsto \phi^*(V)$$ along a group homomorphism $\phi:G\to F$ is tensor.
Clearly $\phi^*\circ\psi^* = (\psi\phi)^*$ for group homomorphisms $\phi:G\to H$ and $\psi:H\to F$. 
In particular group automorphisms of $G$ give rise to tensor autoequivalences of $\V(G)$.
Thus we have a homomorphism of groups
\beq\lb{ags}Aut(G)\to Aut_\ot(\V(G))\qquad \phi\mapsto (\phi^{-1})^*\ \eeq 
Clearly the effect of $\phi^*$ on isomorphism classes of simple objects is $\phi^{-1}$:
$$\phi^*(I(f)) = I(\phi^{-1}(f))\ .$$
Thus the homomorphism \eqref{ags} is a splitting for the epimorphism \eqref{aos}.
In other words $Aut_\ot(\V(G))$ is the semi-direct product of $Aut(G)$ and the kernel of \eqref{aos}.

Now we describe the kernel of \eqref{aos}. 
By the remark \ref{fuc} the kernel coincides with the group  $Aut^1_\ot(\V(G))$ of soft tensor autoequivalences of $\V(G)$.
By the remark \ref{twf} a soft tensor autoequivalence of $\V(G)$ is a twisted form of the identity functor.
Let $\gamma:G^{\times 2}\to k^*$ be a function. Here $k^* = k\setminus\{0\}$ is the multiplicative group of the ground field $k$. 
For $U,V\in \V(G)$ define an isomorphism
\beq\lb{tms}Id(\gamma)_{U,V}:U\ot V\to U\ot V,\qquad u\ot v\mapsto \gamma(f,g)(u\ot v),\quad u\in U_f,\ v\in V_g\ .\eeq
Clearly the collection $Id(\gamma)_{U,V}$ is natural in $U,V\in \V(G)$. Moreover it is the most general form of a natural collection of isomorphisms $U\ot V\to U\ot V$.
\nl
Consider $Id(\gamma)_{U,V}$ as the tensor structure for the identity functor $Id:\V(G)\to \V(G)$. 
The coherence axiom for the tensor structure of $Id(\gamma)_{U,V}$ is equivalent to the equation 
$$\gamma(f,gh)\gamma(g,h) = \gamma(fg,h)\gamma(f,g)$$ for any $f,g,h\in G$. 
Thus $\gamma$ is a 2-cocycle of $G$ with coefficients in $k^*$. 
\nl
Unit axioms are equivalent to the normalisation condition:
$$\gamma(f,e)=\gamma(e,g)=1\qquad f,g\in G\ .$$
\nl
Thus $Id(\gamma)$ with $\gamma\in Z^2(G,k^*)$ is a (general) twisted form of the identity functor.
The composition of twisted forms corresponds to the product of cocycles:
$$Id(\gamma)\circ Id(\gamma') = Id(\gamma\gamma')\ .$$
\nl
A natural transformation $Id\to Id$ amounts to a function $c:G\to k^*$: 
$$c_V:V\to V,\qquad v\mapsto c(g)v,\quad v\in V_g\ .$$
A natural transformation $Id(\gamma)\to Id(\gamma')$ is tensor if the function $c:G\to k^*$ is a coboundary for $\gamma(\gamma')^{-1}$:
$$\gamma(f,g)c(fg) = c(f)c(g)\gamma'(f,g)\ .$$
Thus we have the following.
\bpr\lb{ags}
The group of isomorphism classes of tensor autoequivalences of the category $\V(G)$ is the semi-direct product:
$$Aut_\ot(\V(G))\ \simeq\ Aut(G)\ltimes H^2(G,k^*)\ .$$
In particular the group of soft tensor autoequivalences of $\V(G)$ is the Schur multiplier of $G$:
$$Aut^1_\ot(\V(G))\ \simeq\  H^2(G,k^*)\ .$$
\epr

\bre\lb{cia}
With a invertible object $P$ in a tensor category $\C$ one can associate a tensor autoequivalence
$$P\ot\ \ot P^{-1}:\C\to\C,\qquad X\mapsto P\ot X\ot P$$
the {\em inner autoequivalence} corresponding to $P\in\C$.
The quotient group of $Aut_\ot(\C)$ by inner autoequivalences is denoted by $Out_\ot(\C)$. 

Clearly the inner tensor autoequivalence $I(g)\ot\ \ot I(g^{-1})$ corresponding to the invertible object $I(g)\in\V(G)$ coincides with the inverse image $\phi_g^*$ along the {\em inner} automorphism $\phi_g(x) = g^{-1}xg$.
This give a categorical action (the {\em adjoint action}) of the group $G$ (the group of isomorphism classes of invertible simple objects) on the category $\V(G)$. The group of outer tensor autoequivalences of $\V(G)$ is
$Out_\ot(\V(G))\ \simeq\ Out(G)\ltimes H^2(G,k^*)\ .$
\ere

\subsection{Soft braided tensor autoequivalences of categories of representations}\lb{sar}

For a finite group $G$ denote by $\Rep(G)$ the category of finite dimensional representations. 
The category $\Rep(G)$ is a braided (symmetric) tensor category. 
Representation categories are contravariant in $G$: for a group homomorphism $\phi:G\to F$ there is a braided tensor functor
$$\phi^*:\Rep(F)\to\Rep(G)$$
called the {\em inverse image} along $\phi$. More precisely for an $F$ representation $V$ the $G$-action on $\phi^*(V)$ has the form
$g(v) = \phi(g)(v)$, where $g\in G$ and $v\in V$. 
Clearly $\phi^*\circ\psi^* = (\psi\phi)^*$ for group homomorphisms $\phi:G\to H$ and $\psi:H\to F$. 
In particular group automorphisms of $G$ give rise to braided tensor autoequivalences of $\Rep(G)$.
Inner automorphisms (automorphisms of the form $\phi(x) = gxg^{-1}$) are tensor isomorphic to the identity functor.
Thus we have a homomorphism of groups
\beq\lb{arc}Out(G)\to Aut_{br}(\Rep(G))\qquad \phi\mapsto (\phi^{-1})^*\ \eeq
Here $Out(G)$ is the group of outer automorphisms of $G$, that is the quotient of the group of automorphisms $Aut(G)$ by its normal subgroup $Inn(G)$ consisting of inner automorphisms. 

\bpr\lb{arc}
The map \eqref{arc} induces an isomorphism
$$Aut_{br}(\Rep(G))\ \simeq\ Out(G)\ .$$
\epr
\bpf
It follows form the Deligne's theorem (on the existence of a fibre functor) \cite{de}.
\epf
\bre\lb{gaa}
The results of \cite{da} provide a more elementary proof.
It was proved in \cite{da} that tensor autoequivalences of $\Rep(G)$ correspond to $G$-biGalois algebras.
It is not hard to see that braided tensor autoequivalences of $\Rep(G)$ correspond to commutative $G$-biGalois algebras.
Being semi-simple these algebras have to be isomorphic to the function algebra $k(G)$ with the left and right $G$-actions given by
$$(ga)(x) = a(xg),\qquad (ag)(x) = a(\phi(g)x),\qquad g,x\in G,\ a\in k(G)\ ,$$
where $\phi:G\to G$ is an isomorphism.
\ere

A group automorphism $\phi:G\to G$ is {\em class-preserving} if it preserves conjugacy classes of $G$, that is for any $x\in G$ there is $g\in G$ such that $\phi(x) = gxg^{-1}$. 
It is straightforward that class-preserving automorphisms are closed under the composition.
Denote by $Aut_{cl}(G)$ the group of class-preserving automorphisms of $G$.
Clearly inner automorphisms are class-preserving. The quotient $Aut_{cl}(G)/Inn(G)$ is denoted $Out_{cl}(G)$. 
\bpr\lb{cpa}
The group of isomorphism classes of soft braided tensor autoequivalences of the category $\Rep(G)$ is isomorphic to the group of outer class-preserving automorphisms of $G$:
$$Aut^1_{br}(\Rep(G))\ \simeq\ Out_{cl}(G)\ .$$
\epr
\bpf
By the remark \ref{bfc} it is enough to prove that 
$Out_{cl}(G)$ is the kernel of the natural action of $Out(G)$ on the Grothendieck ring $K_0(\Rep(G))$.
It is well known (from the character theory) that $K_0(\Rep(G))$ is naturally isomorphic to the character ring $R(G)$.
Moreover the character ring $R(G)$ is a subring of the algebra $R(G)\ot_\bZ k$ and the algebra $R(G)\ot_\bZ k$ coincides with the algebra of $k$-valued class functions (functions constant on conjugacy classes of $G$).
Clearly the kernel of the natural action of $Out(G)$ on class functions is $Out_{cl}(G)$. 
\epf

\subsection{Soft braided tensor autoequivalences of Drinfeld centres of finite groups}\lb{sbadc}

We call a $G$-action on a $G$-graded vector space $V = \oplus_{g\in G}V_g$ {\em compatible} (with the grading) if $f(V_g) = V_{fgf^{-1}}$.
Let $\Z(G)$ be the category $G$-graded vector spaces with compatible $G$-actions and with morphism being linear maps preserving grading and action.
Define the tensor product $V\ot U$ of objects $V,U\in\Z(G)$ as the tensor product of $G$-graded vector spaces with the charge conjugation $G$-action.
The category $\Z(G)$ is a tensor category with respect to this tensor product and the trivial associativity constraint.
Moreover $\Z(G)$ is braided with the braiding
$$c_{V,U}(v\otimes u) = f(v)\otimes u,\qquad v\in V_f,\ u\in U\ .$$
We call the category $\Z(G)$ the {\em Drinfeld centre} of $G$ (the author first learned about Drinfeld's work on these categories from \cite{le}).

The functor 
$$\Rep(G)\to \Z(G)$$
considering a $G$-representation as the trivially $G$-graded (concentrated in the trivial degree) is a braided tensor fully faithful functor (full embedding).
This functor has a tensor splitting, more precisely the functor
$$\Z(G)\to \Rep(G)$$ forgetting $G$-grading is tensor and the composition $\Rep(G)\to \Z(G)\to \Rep(G)$ is the identity.
The functor forgetting $G$-action
$$\Z(G)\to \V(G)$$ is also tensor.

Recall that the monoidal centre $\Z(\V(G))$ of the category of graded vector spaces is braided equivalent to $\Z(G)$.
Indeed an object $Z\in \Z(\V(G))$ has a natural half-braiding in $U\in \V(G)$:
$$z_U:U\ot Z\to Z\ot U,\qquad z_U(u\ot z) = f(z)\ot u,\qquad u\in U_f,\ z\in Z$$
and the braided tensor functor 
$$\Z(G)\to\Z(\V(G)),\qquad Z\mapsto (Z,z)$$ is an equivalence (see e.g. \cite{da} for details).
By the functoriality of the monoidal centre there is a homomorphism
\beq\lb{gaz}Aut_\ot(\V(G))\to Aut_{br}(\Z(G))\qquad F\mapsto \tilde F\ .\eeq
Here $\tilde F:\Z(G)\to\Z(G)$ is the braided tensor functor associated with a tensor functor $F:\V(G)\to\V(G)$ as follows
$\tilde F(Z,z) = (F(Z),F(z))$, where $F(z)_U$ is given by
$$\xymatrix{F(U)\ot F(Z) \ar[r]^{F_{U,Z}} & F(U\ot Z) \ar[r]^{F(z_U)} & F(Z\ot U) \ar[r]^{F_{Z,U}^{-1}} & F(Z)\ot F(U)}$$
The tensor structure of $\tilde F$ is the tensor structure of $F$: $\tilde F_{(Z,z),(Z',z')} = F_{Z,Z'}$. 

Thus there is a homomorphism of groups
\beq\lb{agc} Aut(G)\ltimes H^2(G,k^*)\simeq Aut_\ot(\V(G)) \to Aut_{br}(\Z(G))\eeq
We can describe braided autoequivalences $\tilde F_{\phi,\gamma}:\Z(G)\to\Z(G)$ explicitly.
First note that for a group isomorphism $\phi:G\to F$ there is a braided tensor equivalence
$$\phi^*:\Z(F)\to\Z(G)$$
called the {\em inverse image} along $\phi$. For $V\in\Z(F)$ the $G$-action on $\phi^*(V)$ has the form
$g(v) = \phi(g)(v)$, where $g\in G$ and $v\in V$ and the $G$-grading is defined by $V_g = V_{\phi^{-1}(g)}$. 
Clearly $\phi^*\circ\psi^* = (\psi\phi)^*$ for group isomorphisms $\phi:G\to H$ and $\psi:H\to F$. 
This gives the homomorphism
\beq\lb{aaz} Aut(G)\ \to\ Aut_{br}(\Z(G))\ .\eeq
It is easy to see that the homomorphism \eqref{aaz} factors through the outer automorphism group $Out(G)$. 

For $\gamma\in Z^2(G,k^*)$ the tensor autoequivalence $F_\gamma:\V(G)\to\V(G)$ induces a braided tensor autoequivalence 
$\tilde F_\gamma:\Z(G)\to\Z(G)$. The autoequivalence $\tilde F_\gamma$ does not change the $G$-grading of $V\in\Z(G)$ but changes the $G$-action to
$$f*v = \frac{\gamma(f,g)}{\gamma(g,f)}f(v),\qquad v\in V_g.$$
The braided autoequivalence $\tilde F_{\phi,\gamma}:\Z(G)\to\Z(G)$ is the composition $F_\gamma\circ\phi^*$.

The monoidal centre $\Z(\Rep(G))$ of the representation category is also braided equivalent to $\Z(G)$.
Again an object $Z\in \Z(\V(G))$ has an (inverse) half-braiding natural in $U\in\Rep(G)$:
$$z_U:Z\ot U\to U\ot Z,\qquad z_U(z\ot u) = f(u)\ot z,\qquad u\in U,\ z\in Z_f\ .$$
Similarly we have a homomorphism
$$Aut_\ot(\Rep(G))\to Aut_{br}(\Z(G))\ .$$
In particular we have a homomorphism $Out(G)\simeq Aut_{br}(\Rep(G))\ \to\ Aut_{br}(\Z(G))$ which coincides with (the factorisation of) \eqref{aaz}.

The following statement was proved in \cite[Corollary 6.9]{nr} and characterises the image of \eqref{agc}.
\bpr\lb{rtr}
The group $Aut_{br}(\Z(G),\Rep(G))$ of the group $Aut_{br}(\Z(G))$ consisting of tensor isomorphism classes of braided tensor autoequivalences of $\Z(G)$ preserving the subcategory $\Rep(G)\to \Z(G)$ is $Out(G)\ltimes H^2(G,k^*)$. 
\nl
That is braided tensor autoequivalence of $\Z(G)$ preserving the subcategory $\Rep(G)\to \Z(G)$ has the form $\tilde F_{\phi,\gamma}$ for a group automorphism $\phi:G\to G$ and a 2-cocycle $\gamma\in Z^2(G,k^*)$. 
\epr

\bre\lb{chz}
Here we briefly recall the character theory of the category $\Z(G)$. We use slightly different notation comparing e.g. to \cite{dpr}.
The {\em character} of an object $V$ of $\Z(G)$ is a function $\chi_V:\{(f,g)\in G^{\times 2}|\ fg=gf\}\to k$ defined by 
$$\chi_V(f,g) = tr_{V_f}(g),$$ where $tr_{V_f}(g)$ is the trace of the linear operator $g:V_f\to V_f,\quad v\mapsto g(v)$.
It is straightforward that a character is a {\em double class function}, that is 
$$\chi(hfh^{-1},hgh^{-1}) = \chi(f,g),\qquad \forall\ h\in G\ .$$
The character of the tensor product $V\ot U$ has the following expression in terms of the characters of $V,U$:
$$\chi_{U\ot V}(f,g) = \sum_{f_1f_2=f}\chi_U(f_1,g)\chi_V(f_2,g)$$
where the sum is taken over all elements $f_1,f_2$ in the centraliser $C_G(g)$.
\ere

Here we have the main result of the paper describing soft braided tensor autoequivalences of $\Z(G)$.
\bth
The group $Aut^1_{br}(\Z(G))$ of tensor isomorphism classes of soft braided tensor autoequivalences of $\Z(G)$ is the subgroup of the semi-direct product $Out(G)\ltimes H^2(G,k^*)$ consisting of such pairs $(\phi,\gamma)$ that
\beq\lb{sba}\chi(\phi(f),\phi(g)) = \frac{\gamma(f,g)}{\gamma(g,f)}\chi(f,g),\qquad f,g\in G\eeq
for all double class functions $\chi$.
\eth
\bpf
By lemma \ref{rsa} a soft braided  tensor autoequivalence $F$ of $\Z(G)$ preserves the subcategory $\Rep(G)\to \Z(G)$.
By proposition \ref{rtr} $F$ is tensor isomorphic to $\tilde F_{\phi,\gamma}$ for a group automorphism $\phi:G\to G$ and a 2-cocycle $\gamma\in Z^2(G,k^*)$. 

By the remark \ref{bfc} $F:\Z(G)\to\Z(G)$ is soft if and only if acts trivially on the Grothendieck ring $K_0(\Z(G))$.
According to the remark \ref{chz} $K_0(\Z(G))$ embeds in the algebra $K_0(\Z(G))\ot_\bZ k$ and the algebra $K_0(\Z(G))\ot_\bZ k$ coincides with the algebra of $k$-valued double class functions.
\nl
We have the following formula for the character of $\tilde F_{\phi,\gamma}(V)$: 
$$\chi_{\tilde F_{\phi,\gamma}(V)}(f,g) = \frac{\gamma(f,g)}{\gamma(g,f)}\chi_V(\phi(f),\phi(g))\ .$$
Thus $\tilde F_{\phi,\gamma}$ is soft if and only if the pair $\phi,\gamma$ satisfies the condition \eqref{sba}.
\epf

We call a group automorphism $\phi:G\to G$ {\em doubly class-preserving} if it preserves conjugacy classes of commuting pairs of elements of $G$, that is for any $x,y\in G$ such that $xy=yx$ there is $g\in G$ such that $\phi(x) = gxg^{-1},\ \phi(y) = gyg^{-1}$. 
It is straightforward that doubly class-preserving automorphisms are closed under the composition.
Denote by $Aut_{2-cl}(G)$ the group of doubly class-preserving automorphisms of $G$.
Clearly inner automorphisms are doubly class-preserving. The quotient $Aut_{2-cl}(G)/Inn(G)$ is denoted $Out_{2-cl}(G)$.

Thus $\tilde F_\gamma:\Z(G)\to\Z(G)$ is soft if and only if 
\beq\lb{dbm}\frac{\gamma(f,g)}{\gamma(g,f)} = 1\qquad\mbox{for any commuting}\ f,g\in G\ .\eeq
Denote by $B(G)$ the subgroup of $H^2(G,k^*)$ consisting of classes of $\gamma$ satisfying \eqref{dbm}.

\bco
The semi-direct product $Out_{2-cl}(G)\ltimes B(G)$ is a subgroup of the group $Aut^1_{br}(\Z(G))$ of tensor isomorphism classes of soft braided tensor autoequivalences of $\Z(G)$.
\eco

\bco
The group $Aut^1_{br}(\Z(G))$ of tensor isomorphism classes of soft braided tensor autoequivalences of $\Z(G)$ fits into an exact sequence
$$\xymatrix{1\ar[r] & B(G) \ar[r] & Aut^1_{br}(\Z(G)) \ar[r] & Out_{cl}(G)}$$
\eco
\bpf
The image of the homomorphism $Aut^1_{br}(\Z(G)) \to Out(G)$ restricting to the subcategory $\Rep(G)\to \Z(G)$ is clearly in the group 
$Aut^1_{br}(\Rep(G))\simeq Out_{cl}(G)$ of soft braided tensor autoequivalences of $\Rep(G)$ which is isomorphic to the group of outer automorphisms of $G$.
The kernel of the restriction homomorphism consists of classes of $\tilde F_\gamma$ with $\gamma\in Z^2(G,k^*)$ satisfying the condition \eqref{dbm} with trivial $\phi$. That is $[\gamma]\in B(G)$. 
\epf

\bco
For a simple $G$ the group $Aut^1_{br}(\Z(G))$ is trivial.
\eco

\section{Examples}\lb{exa}

Here we give examples of finite groups $G$ with non-trivial groups of soft braided autoequivalences $Aut^1_{br}(\Z(G))$. 

\subsection{Bogomolov multipliers}

It is well-known that for an abelian $A$ the map
$$H^2(A,k^*)\to Hom(\Lambda^2A,k^*),\qquad \gamma\mapsto \Big(\ a\wedge b\mapsto \frac{\gamma(a,b)}{\gamma(b,a)}\ \Big)$$ is an isomorphism.
Thus the group $B(G)$ can be described as the kernel
$$ker\big(H^2(G,k^*) \stackrel{\small res}{\longrightarrow}  \bigoplus_{A\subset G} H^2(A,k^*)\big)$$
of restriction homomorphisms, where the sum is taken over all (2-generated) abelian subgroups $A\subset G$.

It was proved in \cite{bog} that $B(G)$ coincides with the unramified Brauer group $H^2_{nr}(k(V)^G,k^*)$, where $V$ is a faithful representation of $G$ over $k$. Non-triviality of this group provides an obstruction to stable rationality of $k(V)^G$. 
The group (denoted $B_0(G)$) was called the {\em Bogomolov multiplier} of $G$ in \cite{ku}. 
Recently quite a few examples of finite $G$ with no-trivial Bogomolov multiplier were presented in the literature.
We copy some here.

Consider the group 
$$G = \Big\langle a,b,c\ |\ a^2=b^2=1,\ c^2=[a,c],\quad [c,b] = [[c,a],a],\quad [[b,a],G] = 1,\quad [G,[G,[G,G]]]\ \Big\rangle\ .$$
of class 3 and order 64. 
It was proved in \cite{mo} that $B(G)\simeq \bZ/2\bZ$.

The following two groups \cite{jm} are metabelian (of class 2) of order $|G| = p^7$ and of period $p$.
$$G = \Big\langle a,b,c,d\ |\ a^p=b^p=c^p=d^p=1,\quad [a,b] = [c,d],\quad [b,d] = [a,b]^\ve[a,c]^\omega,\quad [G,[G,G]]\ \Big\rangle\ ,$$
where $\ve=1$ for $p=2$ and $\ve=0$ for odd primes $p$, and $\omega$ is a generator of the multiplicative group $(\bZ/p\bZ)^*$.
\nl
The Bogomolov multiplier is $B(G)\simeq \bZ/p\bZ\times \bZ/p\bZ$.
$$\Big\langle a,b,c,d\ |\ a^p=b^p=c^p=d^p=1,\quad [a,b] = [c,d],\quad [a,c] = [a,d] = 1,\quad [G,[G,G]]\ \Big\rangle\ ,$$
\nl
The Bogomolov multiplier is $B(G)\simeq \bZ/p\bZ$.

\subsection{Double class-preserving automorphisms}

First examples of non-inner double class-preserving automorphisms were given in \cite{ne}.
Here we reproduce the construction along with some examples from \cite{Sz}.

Let $M$ be a (finite) abelian group. Consider the action of the additive group $End(M)$ on $M\op M$:
$$f(x,y) = (x,f(x)+y),\qquad f\in End(M),\ x,y\in M\ .$$
Let $G=End(M)\ltimes (M\op M)$ be the semi-direct product with respect to this action.
Note that 
$$[G,G] = Z(G) = 0\ltimes (0\op M)\ .$$
In particular the group $G$ is metabelian.
\nl
For an additive subgroup $E\subset End(M)$ consider the (normal) subgroup
$$G(E) = E\ltimes (M\op M)\ \subset\ G\ .$$
Consider
$$\tilde E = \{g\in End(M)|\ \forall x,y\in M\ \exists f\in E:\ f(x)=g(x),\ f(y)=g(y)\}\ .$$
It was argued in \cite{Sz} that the quotient $\tilde E/E$ embeds into $Out_{2-cl}(G(E))$. 
\nl
Now following \cite{Sz} take $M$ to be an $n$-dimensional vector space over a finite field $\bF_q$ and take 
$E = \fsl(M) = \{g\in End(M)|\ tr(g)=0\}$. 
It follows from double transitivity of the $\fsl(M)$-action on $M$ for $n\geq 3$ that $\widetilde{\fsl(M)}=End(M)$.
\nl
Thus $Out_{2-cl}(G(\fsl(M)))$ has an element of order $q$.

\section{Applications}\lb{app}

Here we sketch a construction of different conformal field theories with the same charge conjugation modular invariant. 

\subsection{Lagrangian algebras of braided autoequivalences}

Let $F:\C\to\C$ be a braided tensor autoequivalence of a braided fusion category $\C$.
Consider an object 
$$Z(F) = \bigoplus_{X\in Irr(\C)}X\boxtimes F(X)^*\in\C\boxtimes\overline\C\ ,$$
where the sum is taken over isomorphism classes of simple objects of $\C$. 
Define a morphism $\mu:Z(F)\ot Z(F)\to Z(F)$ in $\C\boxtimes\overline\C$ as the image of (the sum of) canonical elements under the map
\begin{align*}
& \bigoplus_{X,Y,Z\in Irr(\C)}\C(X\ot Y,Z)\ot_k\C(X\ot Y,Z)^* \simeq \\
& \bigoplus_{X,Y,Z\in Irr(\C)}\C(X\ot Y,Z)\ot_k\C(Z,X\ot Y) \simeq \\
& \bigoplus_{X,Y,Z\in Irr(\C)}\C(X\ot Y,Z)\ot_k\overline\C(F(Z),F(X\ot Y)) \simeq \\
& \bigoplus_{X,Y,Z\in Irr(\C)}\C(X\ot Y,Z)\ot_k\overline\C(F(Z),F(Y)\ot F(X)) = \\
& \bigoplus_{X,Y,Z\in Irr(\C)}\C(X\ot Y,Z)\ot_k\overline\C(F(X)^*\ot F(Y)^*,F(Z)^*) = \\
& \bigoplus_{X,Y,Z\in Irr(\C)}(\C\boxtimes\overline\C)\big((X\ot Y)\boxtimes (F(X)\ot F(Y)),Z\boxtimes F(Z)\big) = \\
& (\C\boxtimes\overline\C)(Z(F)\ot Z(F),Z(F))
\end{align*}
The first isomorphism comes from the non-degenerate pairing 
$$\C(X,W)\ot_k\C(W,X)\to\C(X,X)\simeq k$$
given by the composition of morphisms in $\C$ (here $X$ is simple);
the second isomorphism is the effect on morphisms of the functor $F$;
the third isomorphism is induced by the inverse to the composition 
$$\xymatrix{F(Y)\ot F(X)\ar[rr]^{F_{Y,X}} && F(Y\ot X) \ar[rr]^{F(c_{Y,X})} && F(X\ot Y)}$$

\ble
The pair $(Z(F),\mu)$ is a commutative algebra in $\C\boxtimes\overline\C$. 
\ele
It follows from the results of \cite[section 3.2]{dno} that $Z(F)$ is a Lagrangian algebra in $\C\boxtimes\overline\C$. 
\nl
Note that the algebras $Z(F), Z(F')$ are isomorphic if and only if the autoequivalences $F, F'$ are tensor isomorphic.

For a Lagrangian algebra $Z\in\C\boxtimes\overline\C$ the class $[Z]$ in the Grothendieck ring $K_0(\C\boxtimes\overline\C)\simeq K_0(\C)\times_\bZ K_0(\C)$ written in the basis of classes of simple objects of $\C$
$$[Z] = \sum_{\chi,\xi}M_{\chi\xi}\chi\boxtimes\overline\xi^*$$
gives rise to a non-negative integer matrix $M=(M_{\chi\xi})$ called the {\em modular invariant} of $Z$.
\nl
We say that $Z$ has the {\em charge conjugation} modular invariant if the matrix $M$ is the identity $M_{\chi\xi} = \delta_{\chi\xi}$.

The next lemma follows from the results of \cite[section 3.2]{dno}.
\ble
A Lagrangian algebra $Z\in\C\boxtimes\overline\C$ has the charge conjugation modular invariant if and only if $Z \simeq Z(F)$ for a soft braided tensor autoequivalence $F:\C\to\C$. 
\ele

\subsection{Holomorphic permutation orbifolds}

Let $V$ be a holomorphic vertex operator algebra (for example $V = \e_{8,1}$). 
\nl
Let $G\subset S_n$ be a subgroup of the permutation group.
The vertex operator subalgebra $(V^{\ot n})^G$ of invariants is called the {\em chiral permutation orbifold} of $V$.
\nl
According to \cite{Ki} its category of representation is $\Rep((V^{\ot n})^G) = \Z(G)$ (subject to the rationality of $(V^{\ot n})^G$). 
\nl
According to \cite{frs4,HK} a Lagrangian algebra $Z\in \Z(G)\boxtimes\overline\Z(G)$ gives rise to a rational conformal field theory with the left (right) chiral algebras $(V^{\ot n})^G$ and the modular invariant $[Z]$.
In particular taking $Z=Z(F)$ for a soft braided tensor autoequivalence $F:\Z(G)\to\Z(G)$ provides an example of non-trivial rational conformal field theory with the charge conjugation modular invariant.

\end{document}